\documentclass[10pt,a4paper]{article} 

\usepackage[latin1,utf8]{inputenc}
\usepackage[T1]{fontenc} 
\usepackage{textcomp} 
\usepackage[english]{babel}
\usepackage[autolanguage]{numprint} 
\usepackage{hyperref} 
\usepackage{graphicx} 
\usepackage{verbatim} 
\usepackage[all]{xy}
\usepackage{float}

\usepackage{lmodern}
\usepackage{mathrsfs}

\usepackage{amsmath,amssymb,amsthm} 
\usepackage{tikz,tkz-tab} 

\usepackage{mathdots}
\usepackage{makeidx}
\usepackage{stmaryrd} 

\usepackage{fancyhdr} 
\renewcommand\epsilon\varespilon 
\renewcommand\phi\varphi 

\newcommand\NN{\mathbb{N}} 
\newcommand\ZZ{\mathbb{Z}} 
\newcommand\QQ{\mathbb{Q}} 
\newcommand\RR{\mathbb{R}} 
\newcommand\CC{\mathbb{C}} 
\newcommand\LL{\mathbb{L}} 
\newcommand\OK{\mathcal{O}_K} 

\newcommand{\wj}[3]{\omega_{#1,#2}(#3)} 
\newcommand{\wjhat}[3]{\widehat{\omega}_{#1,#2}(#3)} 

\newcommand\bx{\mathbf{x}}  
\newcommand\by{\mathbf{y}}  
\newcommand\bu{\mathbf{u}}  

\newcommand{\psc}[2]{\left\langle#1,#2\right\rangle}

\newcommand{\Vect}[2][]{\textrm{ {\rm Span}}_{#1}\left(#2\right)} 
\newcommand{\dist}[2]{\textrm{ {\rm dist}}\left(#1,#2\right)}

\newcommand{\distI}[3]{\omega_{#1}\left(#2,#3\right)}

\theoremstyle{definition} 
\newtheorem{Def}{Definition}

\theoremstyle{plain} 
\newtheorem{Prop}[Def]{Proposition} 
\newtheorem{Thm}[Def]{Theorem} 
\newtheorem{Thm2}{Theorem}

\theoremstyle{remark} 
\newtheorem{Rem}[Def]{Remark} 

\newenvironment*{Dem}{\noindent{\bf Proof}}{\hfill$\square$\bigskip$\newline$}

\numberwithin{equation}{section} 

\newcommand{\Addresses}{{
  \bigskip
  \footnotesize

  A.~Poëls, \textsc{Laboratoire de Math\'ematiques d'Orsay, Univ. Paris-Sud, CNRS,
Universit\'e Paris-Saclay, 91405 Orsay, France}\par\nopagebreak
  \textit{E-mail} : \texttt{anthony.poels@math.u-psud.fr}

}}

\newcommand{\MSC}{{
  \footnotesize

  \textbf{MSC~2010}: 11K60 (Primary); 11J99 (Secondary).
}}

\newcommand{\keysW}{{
  \footnotesize

  \textbf{Keywords}: transference theorems; diophantine approximation exponents; approximation of subspaces
}}

\newcommand{\Ack}{{
  \footnotesize

  \textbf{Acknowledgements}: I would like to thank Yann Bugeaud, Michel Laurent, Damien Roy, Michel Waldschmidt
  and the referee for their useful comments. I am also very grateful to Stéphane Fischler for giving me a lot of
  feedback on this work.\\
  
  The final publication is available at Springer via http://dx.doi.org/10.1007/s00605-017-1084-4
}}

\title{The complex case of Schmidt's going-down Theorem}  
\author{Anthony Poels} 

\begin{document} 
\maketitle

\begin{abstract}
  In  $~1967$, Schmidt wrote a seminal paper \cite{Schmidt67} on heights of subspaces of $\RR^n$ or $\CC^n$ defined
  over a number field $K$, and diophantine approximation problems. The going-down Theorem -- one of the main
  theorems he proved in his paper -- remains valid in two cases depending on whether the embedding of $K$ in the
  complex field $\CC$ is a real or a complex non-real embedding. For the latter, and more generally as soon as $K$ is not
  totally real, at some point of the proof, the
  arguments in \cite{Schmidt67} do not exactly work as announced. In this note, Schmidt’s ideas are worked out in
  details and his proof of the complex case is presented, solving the aforementioned problem. Some definitions of
  Schmidt are reformulated in terms of multilinear algebra and wedge product, following the approaches of
  Laurent \cite{laurent703146transfer}, Bugeaud and Laurent \cite{bugeaud2010transfer} and Roy
  \cite{Roy2014ConstructionRegularSystem}, \cite{Roy_octobre}.\\
  In \cite{laurent703146transfer}
  Laurent introduces in the case $K=\QQ$ a family of exponents and he gives a series of inequalities relating them.
  In Section \ref{section application exponents} these exponents are defined for an arbitrary
  number field $K$. Using the going-up and the going-down Theorems Laurent's inequalities are generalized to this setting.
\end{abstract}

\MSC

\keysW

\Ack

\section{Introduction}

  In a paper \cite{Schmidt67} written in $1967$, Schmidt generalizes the basic diophantine approximation
  problem « given a real number $\alpha$, how « well » can it be approximated by rational numbers ? » as follows.
  Let $A$ be a subspace of a Euclidean or unitary space $G^n$ of dimension $n$. Suppose that $A$ has dimension
  $0 < d < n$. How « well » can $A$ be approximated by subspaces $B$ of dimension $e$ defined over a given number
  field $K$ ? Formulating precisely what « well » means requires some work. Schmidt binds two different notions
  that are recalled in Section \ref{section preliminaires} below : $A$ is « well » approximated by $B$ if on the
  one hand $A$ and $B$ are « close » (Schmidt uses several angles of inclination to measure this « closeness »,
  cf. Proposition \ref{prop def des lambda i} and following definitions), and on the other hand $B$ is not too
  « complicated » (Schmidt uses the notion of the \emph{height} of a subspace to measure its « complicatedness »,
  cf. \eqref{Def premiere hauteur} and \eqref{Prop deuxieme def hauteur}).\\
  In his article, Schmidt establishes several transference theorems of the Perron-Khintchine-type (see for
  example \cite{khintchine1925zwei}, \cite{khintchine1926klasse}, \cite{perron1921diophantische}). These theorems
  lead to the conclusion that if  a subspace $A$ can be well approximated by subspaces of dimension $e$ ($0<e<n$),
  then it can also be well approximated by subspaces of any given dimension $e'$. Schmidt's going-down Theorem
  (\cite{Schmidt67} Theorem $10$) is one of these transference theorems (treating the case $e' < e$) and is useful
  to prove diophantine approximation theorems (as \cite{Schmidt67} Theorem $13$ for example). More recently this
  work was revisited by Laurent \cite{laurent703146transfer} and Bugeaud and Laurent \cite{bugeaud2010transfer} in
  the case where $A$ is a one-dimensional subspace of $\RR^n$ and $K = \QQ$. Laurent introduces a family of
  approximation exponents to points in $\RR^n$ by linear subspaces and using going-up and going-down Theorems he
  proves a series of inequalities relating these exponents \cite{laurent703146transfer}. Roy shows \cite{Roy_octobre}
  that the going-up and going-down transference inequalities of Schmidt and Laurent describe the full spectrum of
  these exponents. In Section \ref{section application exponents} these exponents are generalized for an arbitrary
  number field $K$. Using the going-up and the going-down Theorems one shows that Laurent's inequalities remain
  valid for the aforementioned generalized exponents.

  \bigskip
  The going-down Theorem remains valid in both the real case $(a)$ and the complex (non-real) case $(b)$ (see below
  for details). At some point of the proof for $(b)$, Schmidt's arguments do not exactly work as announced; as the
  referee pointed out, this happens also in case $(a)$ if $K$ is not totally real (see  Remarks
  \ref{Remarque preuve Schmidt} and \ref{Remarque2 preuve Schmidt} for more technical details). The main goal of this note is to work out Schmidt's
  ideas in details  and solve this problem.\\
In case $(a)$, $G^n$ denotes Euclidean space $(\RR^n,\psc{\cdot}{\cdot})$, $K$ a number field embedded in $\RR$, and $q=1$.\\
In case $(b)$, $G^n$ denotes unitary space $(\CC^n,\psc{\cdot}{\cdot})$, $K$ a number field embedded in $\CC$ (but this embedding is not real), and $q=2$.\\
In this paper $C_1, C_2,\dots$ will be positive constants depending only on $K$ and $n$ but independent from the
subspaces $A^d, B^e,\dots$ considered. We also keep as far as possible the same numbering as Schmidt (so that the
numbering does not start with $C_1$). The notation $A^d$ (or $B^d,S^d,\dots$) will always mean that
$A^d$ (or $B^d,S^d,\dots$) has dimension $d$.

\begin{Thm2}[Schmidt's going-down Theorem]$\newline$
\label{thm going-down}
Let $A^d$, $B^e$ be subspaces of $G^n$ with $B^e$ defined over $K$ and of height $H(B^e) \leq H$ (where $H\geq 1$ is a fixed constant). Let $1\leq h\leq f'=\min(d,e-1)$, $c\geq 1$ and assume that
\begin{align}
\label{going_down 1ère equation}
H(B^e)\omega_i^q(A^d,B^e)\leq c^qH^{-(qy_i-1)} \quad (i=1,\dots,h).
\end{align}
where $y_1\geq\dots\geq y_h\geq (qh)^{-1}$. Put $y=y_1+\dots+y_h$ and assume
\begin{align}
\label{going_down 2ème equation}
y_i':= y_ie(qy+e-1)^{-1} \geq q^{-1}\quad (i=1,\dots,h).
\end{align}
Then there is a subspace $B^{e-1}\subset B^e$, defined over $K$, of height
\[
H(B^{e-1})\leq C_5H(B^e)H^{(qy-1)/e}\leq C_5H^{(e+qy-1)/e}=:H'
\]
having
\[
H(B^{e-1})\omega_i^q(A^d,B^{e-1})\leq C_6c^qH^{-(qy'_i-1)(qy+e-1)/e}= C_7c^qH'^{-(qy_i'-1)}\quad(i=1,\dots,h),
\]
whence
\[
\omega_i(A^d,B^{e-1})\leq C_8cH(B^{e-1})^{-y_i'}\quad(i=1,\dots,h).
\]
On the other hand, if instead of \eqref{going_down 1ère equation},
\begin{align}
\label{going_down 3ème equation}
\omega_i(A^d,B^e)=0\quad (i=1,\dots,h),
\end{align}
put
\[
y_0':= e(qh)^{-1}.
\]
Then, for any given $H'\geq C_9H$ there is a subspace $B^{e-1}\subset B^e$, defined over $K$, of height $H(B^{e-1})\leq H'$, having
\[
H(B^{e-1})\omega_i^q(A^d,B^{e-1})\leq C_{10}H^{qy_0'}H'^{-(qy_0'-1)}\quad (i=1,\dots,h),
\]
whence
\[
\omega_i(A^d,B^{e-1})\leq C_{11}H^{y_0'}H(B^{e-1})^{-y_0'}\quad (i=1,\dots,h).
\]
In this theorem, $C_5,\dots,C_{11}$ are positive constants which depend on $K, n, y_1,\dots,y_h$ but not on $A^d, B^e, B^{e-1}, H, c$.
\end{Thm2}
In Section \ref{section preliminaires} we recall the definitions of the \emph{height} of a subspace and of the functions
$\omega_i$ which are used to measure the closeness of two subspaces. We also recall results of  \cite{Schmidt67}
which we will need in the proof of Theorem \ref{thm going-down}. Some definitions of Schmidt are reformulated in
terms of multilinear algebra and wedge product -- this is the case of the definition of the height of a subspace for
instance -- following the approaches of Laurent \cite{laurent703146transfer},
Bugeaud and Laurent \cite{bugeaud2010transfer} and Roy \cite{Roy2014ConstructionRegularSystem}, \cite{Roy_octobre}.
In Section \ref{section More notations} we introduce
specific notation for the complex case in order to avoid confusion between the inner product
$\psc{\bx}{\by} = \sum x_i\overline{y_i}$ and the bilinear form $\phi(\bx,\by) = \sum x_iy_i$ (which are denoted in
the same way in \cite{Schmidt67}). In Section \ref{section proof going-down theorem} we present Schmidt's proof of the
going-down Theorem in the complex case, solving in Section \ref{subsection Proporties of W} the problem alluded to above
(see also Remarks \ref{Remarque preuve Schmidt} and \ref{Remarque2 preuve Schmidt} in Section \ref{section de la remark corrective}).
In Section \ref{section application exponents} the generalization of Laurent's exponents is given.

\section{Multilinear algebra, distance and height of subspaces}
  \label{section preliminaires}
  \label{subsection multilinear height distance}

  In this section one reformulates some definitions of Schmidt \cite{Schmidt67} -- among others the height $H(S)$ of a subspace $S$
  and the quantities $\distI{i}{A}{B}$ which characterize the distance between two subspaces $A$ and $B$ --
  in terms of multilinear algebra and wedge product. This approach has already been investigated by Laurent in
  \cite{laurent703146transfer} and Bugeaud and Laurent \cite{bugeaud2010transfer} in order to give another proof
  of the going-up and going-down transfers in the case $K=\QQ$ with $A^d$ of dimension $d=1$.
  See also Roy \cite{Roy2014ConstructionRegularSystem} and \cite{Roy_octobre} for further examples of the use of such
  tools in the context of parametric geometry of numbers.\\

Let $\LL = \RR$ or $\CC$ and let $n\in \NN^*$. We endow $G^n := \LL^n$ with its usual structure of
inner product space. Let $\psc{\cdot}{\cdot}$ denotes the canonical inner product on $G^n$ (if $\LL = \CC$,
we ask for the linearity of the first argument) and $||\cdot||$ its associated norm. If we fix an integer
$m$, $1\leq m \leq n$, we always endow the vector space $\bigwedge^m(G^n)$ with the unique structure of
inner product space such that, for any orthonormal basis $(e_1,\dots,e_n)$ of $G^n$, the products
$e_{i_1}\wedge\dots\wedge e_{i_m}$ ($i_1<\dots <i_m$) form an orthonormal basis of $\bigwedge^m(G^n)$.
We still denote by $\psc{\cdot}{\cdot}$ its inner product and by $||\cdot||$ the associated norm.
(Note that with this notation, we have
$D(X_1,\dots,X_m) = ||X_1\wedge\dots\wedge X_m||$ for any $X_1,\dots,X_m\in G^n$, where
$D(X_1,\dots,X_m) = \Big(\det\Big(\psc{X_i}{X_j}\Big)_{i,j}\Big)^{1/2}$
denotes the \emph{generalized determinant} of $(X_1,\dots,X_m)$, see \cite{Schmidt67} for more details about
this notion). \\

Let $K$ be an algebraic number field of degree $[K:\QQ]=p$ and $\OK$ be the ring of integers of $K$.
Let $\sigma_1,\dots,\sigma_p$ be the different embeddings of $K$ into the field $\CC$ of complex numbers.
For $\xi\in K$, put $\xi^{(i)}$ for the image of $\xi$ under $\sigma_i$. Similarly, if
$X = (\xi_1,\dots,\xi_n) \in K^n$, put $X^{(i)} = (\xi_1^{(i)},\dots, \xi_n^{(i)})$. Let $S^d$ a subspace
of $K^n$ of dimension $d$. Let $(X_1,\dots,X_d)$ be a basis of $S^d$ and form the matrix $M$  with row
vectors $X_1,\dots,X_d$. Let $\mathfrak{a}$ be the fractional ideal of $K$ generated by the $\binom{n}{d}$
determinants of all $d\times d$-- submatrices of $M$. The \emph{height} of $S^d$ is defined by
\begin{equation}
\label{Def premiere hauteur}
    H(S^d) = N(\mathfrak{a})^{-1} \prod_{j=1}^{p}||X_1^{(j)}\wedge\dots\wedge X_d^{(j)}||,
\end{equation}
where $N(\mathfrak{a})=N_{K/\QQ}(\mathfrak{a})\in\QQ^+$ denotes the \emph{norm} of the ideal $\mathfrak{a}$.
This definition does not depend of the choice of the basis $(X_1,\dots,X_d)$. See \cite{Schmidt67} §$1$ for more
explanations about this notion.

We suppose now that $K$ is embedded in $\LL$. We denote by $\bigwedge^m(\OK^n)$ ($1\leq m\leq n$) the free $\OK$-module of rank
$\binom{n}{m}$ spanned by the products $x_1\wedge\dots\wedge x_m$ with $x_1,\dots,x_m\in \OK^n$.
 A subspace $S$ of $\CC^n$ is said to be defined over $K$ if it is defined by
linear equations with coefficients in $K$ (or equivalently, if there is a basis of $S$ with coordinates in $K$).
If $S$ is defined over $K$, one can considerer its height as the height of $S\cap K^n$. If $(X_1,\dots,X_d)\in K^n$ form
a basis of $S^d$, the associated fractional ideal $\mathfrak{a}$ defined above is the fractional ideal generated by the
coordinates of $X_1\wedge\dots\wedge X_d$ with respect to a basis of $\bigwedge^m(\OK^n)$.
Note that if $K = \QQ$ and $\OK = \ZZ$, we may suppose that $(X_1,\dots,X_d)$ form
a basis of $S^d\cap \ZZ^n$ (and so, that it can be extended to a basis $(X_1,\dots, X_n)$ of $\ZZ^n$, which is
equivalent to asking that $X_1\wedge\dots X_d$ is a primitive vector, \emph{i.e} $N(\mathfrak{a}) = 1$).
In the general case since $\OK$ is not necessarily a principal ring, $S^d\cap\OK^n$ may not be a free $\OK$-module.
However, the ideal class group of $K$ is finite and using a system of representatives consisting of integral ideals, it can be proved that
$X_1,\dots,X_d$ may be chosen such that $X_1,\dots,X_d\in\OK^n$ and $N(\mathfrak{a}) \leq C$ where $C>0$ depends of $K$ only.\\
Formula \eqref{Prop deuxieme def hauteur} in §~\ref{section preliminaires} allows one to consider the height of some subspace more geometrically.\\

Finally, one has to introduce the functions $\omega_i$ used by Schmidt to measure the "closeness" of two subspaces $A$ and $B$
of a Euclidean or unitary space. We define the (projective) distance between two non-zero vectors $X$ and $Y$ of $G^n$ by
\[
    \dist{X}{Y} := \frac{||X\wedge Y||}{||X||\; ||Y||}.
\]
Note that in \cite{Schmidt67} $\dist{X}{Y}$ is denoted by $\omega(X,Y)$. It satisfies the triangle inequality
\[
    \dist{X}{Z} \leq \dist{X}{Y} + \dist{Y}{Z}\quad X,Y,Z\in G^n\setminus\{0\}.
\]
For $X\in G^n\setminus\{0\}$ and a subspace $B^e\neq \{0\}$ of $G^n$, we define the distance from $X$ to $B^e$ by
\[
    \dist{X}{B^e} = \inf_{Y\in B^e\setminus\{0\}}\dist{X}{Y}.
\]
Note that
\begin{equation}
\label{Eq omega de l orthogonal}
    \dist{X}{B}^2 + \dist{X}{B^{\perp}}^2 = 1
\end{equation}
for every subspace $B$ of dimension $0 < e < n$ (see \cite[Section~8]{Schmidt67} formula $(8)$).
Note once again that in \cite{Schmidt67} $\dist{X}{B^e}$ is denoted by $\omega(X,B^e)$ and that this infimum is in fact a minimum.

\begin{Def}
Let $A^d$ and $B^e$ be subspaces of $G^n$ of dimensions $d$ and $e$ respectively, with $f:=\min(d,e) > 0$. Set
\[
    \distI{i}{A^d}{B^e} := \inf_{F^i\subset A^d}\sup_{X\in F^i\setminus\{0\}} \dist{X}{B^e},
\]
for $i=1,\dots,f$. Here $F^i$ refers to an arbitrary subspace of $A$ of dimension $i$.
\end{Def}
Intuitively, the smaller the $\omega_i$ are, the closer $A^d$ and $B^e$ are.
These quantities are the same as those introduced by Schmidt in \cite[§8]{Schmidt67} (it is a direct consequence of
Schmidt's definitions, his Lemma~12 and Lagrange's identity \mbox{$||X||^2||Y||^2 = |\psc{X}{Y}|^2 + ||X\wedge Y||^2$}).
In particular, we have the useful following result (see \cite[Theorem $4$ on page $443$]{Schmidt67} noting that
$\lambda_i = \sqrt{1-\omega_i^2}$) :
\begin{Prop}
\label{prop def des lambda i}
Let $A^d$ and $B^e$ be subspaces of $G^n$ of dimensions $d$ and $e$ respectively, with $f:=\min(d,e)>0$. Then there are
orthonormal bases $X_1,\dots,X_d$ and $Y_1,\dots,Y_e$ of $A^d, B^e$ respectively, and reals
$0\leq\omega_1\leq\dots\leq\omega_f\leq 1$ such that
\[
\dist{X_i}{Y_j}=\left\{
\begin{array}{ll}
\omega_i \qquad \textrm{if $i=j$}\\
1 \qquad \textrm{otherwise}
\end{array}
\right.
\qquad (1\leq i\leq d,\;1\leq j\leq e).
\]
The numbers $\omega_1,\dots,\omega_f$ are independent of any freedom of choice in $X_i,Y_j$ and are invariant under
unitary transformations applied simultaneously to $A^d, B^e$. Moreover, one has
\[
    \omega_i = \distI{i}{A^d}{B^e} \quad (1\leq i\leq f).
\]
\end{Prop}

If $d+e\leq n$, set
\[
    \mu(A^d,B^e):=\prod_{k=1}^f\omega_k(A^d,B^e).
\]
(Although we will not use it in this paper, if $d+e > n$, $\mu(A^d,B^e)$ can be defined as $\prod_{k=1}^f\omega_{k+g}(A^d,B^e)$ where
 $g:= d+e-n$. See Sections $7$ and $8$ of \cite{Schmidt67} for more details).

The next and last proposition is an equivalent definition of $\mu$ in the case $d+e\leq n$ (see \cite[§ $6$-$8$]{Schmidt67}, especially formula $(7)$ ).
\begin{Prop}
\label{Prop lien mu et Det generalise}
If $(X_1,\dots,X_d)$ and $(Y_1,\dots,Y_e)$ are arbitrary bases of $A^d, B^e$, respectively, and $d+e\leq n$, one has
\[
    \mu(A^d,B^e)=\frac{||X_1\wedge\dots\wedge X_d\wedge Y_1\wedge\dots\wedge Y_e||}
    {||X_1\wedge\dots\wedge X_d||\;||Y_1\wedge\dots\wedge Y_e||}.
\]
\end{Prop}
This formula generalizes the formula $(4.1)$ of \cite{bugeaud2010transfer}, which describes the special case $d=f=1$
(in this case, $\mu(A^d,B^e) = \dist{X_1}{B^e}$).\\

For $n\in\NN^*$, $E^n$ denotes the Euclidean space $\RR^n$ with its canonical scalar product.
A \emph{lattice} of $E^n$ will mean a discrete group of vectors of $E^n$ (not necessarily cocompact).
The \emph{rank} of a lattice is the maximal number of linearly independent vectors of the lattice. Define the
determinant of a lattice $\Lambda$ of rank $m$ by $d(\Lambda) =  ||X_1\wedge\dots\wedge X_m||$ where $X_1,\dots,X_m$ are
basis vectors of $\Lambda$ if $m>0$, and by $d(\Lambda)= 1$ if $\Lambda=\{0\}$.\\
Suppose now that $K\subset\CC$ but $K\not\subset\RR$. Let $p=r_1+2r_2$ and $\xi^{(2r_2+j)}$ be real for
$1\leq j\leq r_1$, $\xi^{(j+1)}$ the complex conjugate of $\xi^{(j)}$ for $1\leq j\leq 2r_2-1$, $j$ odd, and
every $\xi\in K$. We may assume $\sigma_1$ to be the identity map, $\sigma_2$ the complex conjugate map (such
that $\sigma_1(\xi), \sigma_2(\xi)$ for $\xi\in\CC$ -- not necessarily in $K$ -- can also be considered). Put
\[
    \Delta = 2^{-r_2}|\delta|^{1/2}
\]
where $\delta$ is the \emph{discriminant} of $K$. Set
\[\xi^{[i]}=\left\{
\begin{array}{ll}
\textrm{Re }\xi^{(i)}\qquad &\qquad\;\textrm{if } 1\leq i\leq 2r_2\textrm{ and $i$ odd},\\
\textrm{Im }\xi^{(i)}\qquad &\qquad\;\textrm{if } 1\leq i\leq 2r_2\textrm{ and $i$ even},\\
\xi^{(i)}\qquad &\qquad\;\textrm{if } 2r_2+1\leq i\leq p.
\end{array}
\right.
\]
Here, $\textrm{Re}$ and $\textrm{Im}$ denote real and imaginary parts. Given
$1\leq i\leq p$ and $X=(\xi_1,\dots,\xi_n)\in K^n$, write $X^{(i)} := (\xi_1^{(i)},\dots, \xi_n^{(i)})$ and
$X^{[i]} := (\xi_1^{[i]},\dots, \xi_n^{[i]})$. For $i=1,2$ $X^{(i)}$ and $X^{[i]}$ are defined for any $X\in\CC^n$.
Then, notice that for all $X\in\CC^n$ one has $X = X^{[1]}-iX^{[2]}$. Let $\rho:K^n\rightarrow E^{np}$ be the
$\QQ$-linear map defined by
\begin{align}
\label{Def map rho}
\rho(X) = (X^{[1]},\dots,X^{[p]})\in E^{np}.
\end{align}
It is the same map $\rho$ as the one defined by Schmidt \cite[p. $435$]{Schmidt67} if one rearranges its coordinates;
this does not change the main property \eqref{Prop deuxieme def hauteur} recalled below.\\
Let $S^d$ be a subspace of $K^n$ of dimension $d$, and $\OK(S^d)$ be the subset of all of $X\in S^d$ whose components
are in $\OK$. Then $\Lambda(S^d):=\rho(\OK(S^d))$ is a lattice in $E^{np}$ of rank $dp$ (see
\cite[§$3$]{Schmidt67}). Moreover Theorem $1$ on page $435$ of \cite{Schmidt67} asserts that\\
\begin{equation}
\label{Prop deuxieme def hauteur}
    H(S^d) = \Delta^{-d}d\big(\Lambda(S^d)\big),
\end{equation}
where $d\big(\Lambda\big)$ denotes the determinant of the lattice $\Lambda$.

\section{Specific notation in the complex case}
\label{section More notations}
In the setting of Section \ref{section preliminaires} let us assume now that $K$ is non-real and $G=\CC^n$. We shall
distinguish carefully the sesquilinear scalar product $\psc{X}{Y}$ on $\CC^n$ and the canonical bilinear form
$\phi(X,Y)$ on $\CC^n$ or $K^n$. By contrast, both are denoted by $XY$ in \cite{Schmidt67}. Of course in the real
case $\psc{X}{Y}$ and $\phi(X,Y)$ coincide.
Let $\psc{\cdot}{\cdot} : \CC^n\times\CC^n\rightarrow\CC$ be the inner product defined by
\[
    \psc{(x_i)_i}{(y_i)_i} = \sum_{i=1}^nx_i\overline{y_i}.
\]
If $W$ is a subspace of $\CC^n$, $W^{\perp}$ denotes the orthogonal complement of $W$ for $\psc{\cdot}{\cdot}$.
Let $K'$ be the complex conjugate field of $K$, $K'=\{\overline{z}\;;\; z\in K\}$. If $W$ is defined over $K$, note
that $W^{\perp}$ is defined over $K'$, and generally not over $K$. If $S$ is a subspace of $\CC^n$ defined over $K$,
it follows easily from the definition of the height that $H(S') =  H(S)$, where $S'$ denotes the set of all
$\overline{z}$ with $z\in S$, which is a subspace defined over $K'$. If $\phi : \CC^n\times\CC^n\rightarrow \CC$
denotes the bilinear form defined by
\[
    \phi((x_i)_i,(y_i)_i) = \sum_{k=1}^nx_iy_i,
\]
and $W^{\phi,\perp}$ denotes its orthogonal complement with respect to $\phi$, one can show that
\[
    H(S^{\phi,\perp}) = H(S)
\]
for all subspaces $S$ of $\CC^n$ defined over $K$ (see \cite{Schmidt67} Eq. $(4)$  on page $433$ and \cite{Hodge}
Theorem $1$ on page $294$, although it is not expressed in the same language). We may in particular deduce from the
last statement this useful result :\\
Let $S$ be a subspace of $\CC^n$ defined over $K$. The subspace $S^{\perp}$ is defined over $K'$ and satisfies
\begin{equation}
\label{Prop hauteur de l orthogonal}
    H(S) = H(S^{\perp}).
\end{equation}

\section{Proof of the going-down Theorem in the complex case}
\label{section proof going-down theorem}

\begin{proof}
In this section one proves Schmidt's going-down Theorem (that is, Theorem \ref{thm going-down} in the introduction),
in case $(b)$. In other words, $K$ is a number field embedded in $\CC$, with $K\not\subset \RR$. We can suppose that
$G^n=\big(\CC^n,\psc{\cdot}{\cdot}\big)$; we keep the notation of Sections \ref{section preliminaires} and
\ref{section More notations}. Let $K'$ be the complex conjugate field of $K$, $K'=\{\overline{z}\;;\; z\in K\}$.\\
Write $p = r_1 + 2r_2$. Notation $\xi^{(i)}$, $\xi^{[i]}$, $\rho$ (cf. \eqref{Def map rho})... will be used with
respect to $K'$ ($\sigma_1,\dots,\sigma_p$ denote the different isomorphisms of $K'$ into the field $\CC$ of complex
numbers etc. Notice that Schmidt does not mention the field $K'$ explicitly).\\
Let $B^e$ be a subspace of $\CC^n$ defined over $K$. Schmidt first assumes that \eqref{going_down 1ère equation} holds.
Let $m=n-e$ and $B^{e,\perp}:=(B^e)^{\perp} = \Vect{Z_1,\dots,Z_m}$ with $Z_i\in K'^n$, where $\Vect{T_1,\dots,T_m}$
denotes the subspace generated by $T_1,\dots,T_m$. We are going to follow Schmidt's idea in order to construct a vector
$W\in K'^n\setminus (B^e)^{\perp}$ such that
\begin{align}
\label{going_down def B^(e-1)}
B^{e-1}:=\Vect{W,Z_1,\dots,Z_m}^{\perp}
\end{align}
(which is defined over $K$) has the required properties.\\
Let $\lambda_i:= \Big(1-\omega_i(A^d,B^e)\Big)^{1/2}$ (for $i=1,\dots,f:=\min(d,e)$) and choose orthonormal bases $(X_1,\dots,X_d)$,
$(Y_1,\dots,Y_e)$ of $A^d$ and $B^e$ respectively having $\psc{X_i}{Y_j}=\delta_{ij}\lambda_i$ (such bases are given
by Proposition \ref{prop def des lambda i}), whence for $i=1,\dots,f$ one has $\dist{X_i}{Y_i}=\omega_i(A^d,B^e)$.
Notice that $\psc{Y_i}{Z_j}=0$ ($i=1,\dots,e$ and $j=1,\dots,m$).\\

The lattice $\Lambda(B^{e,\perp})\subset E^{pn}$ (constructed from $K'$ and $\rho$ as in §$2$) has rank $pm$ and
determinant $\Delta^m H(B^{e,\perp})=\Delta^m H(B^e)$ (by \eqref{Prop hauteur de l orthogonal}).
Let $\mathfrak{I}_1,\dots,\mathfrak{I}_{pm}$ be a basis of this lattice and let $\Pi$ be the set of points
$\mathfrak{I} = \sum c_i\mathfrak{I}_i$ with $c_i\in [-1/2,1/2]$. The set $\Pi$ has $pm$-dimensional volume
$\Delta^m H(B^e)$ and contains no lattice point of $\Lambda:=\Lambda(K'^n)$ but $0$. Set
$S^*:=\Vect{\Lambda(B^{e,\perp})}\subset E^{pn}$; it has dimension $pm$.\\\\

\subsection{Construction of $W$}

Remember that $\rho$ is defined with respect to $K'$ (not $K$).
For all $Y\in \CC^n$ and $Z\in K'^n$, one has
\begin{align*}
    \psc{Y}{Z} = \psc{(Y^{[1]},Y^{[2]},0,\dots,0)}{\rho(Z)}+i\psc{(-Y^{[2]},Y^{[1]},0,\dots,0)}{\rho(Z)}.
\end{align*}
$(Y_j)_j$ is an orthonormal basis and each $Y_j$ is orthogonal to $B^{e,\perp}$, hence
\[
    \mathfrak{Y}_j^1:=(Y_j^{[1]},Y_j^{[2]},0,\dots,0) \quad (j=1,\dots,h),
\]
and
\[
    \mathfrak{Y}_j^2:=(-Y_j^{[2]},Y_j^{[1]},0,\dots,0) \quad (j=1,\dots,h),
\]
form an orthonormal family of $2h$ vectors of $E^{pn}$ which is orthogonal to $S^*$. Set
$T_h:=\Vect{\mathfrak{Y}_1^1,\mathfrak{Y}_1^2,\dots,\mathfrak{Y}_h^1,\mathfrak{Y}_h^2}$; then $T_h$ is a
subspace of dimension $2h$. A vector $\mathfrak{X}\in E^{pn}$ can be uniquely written as
\[
    \mathfrak{X} = \mathfrak{X}^*+\mathfrak{X}_T+\mathfrak{X}_0,
\]
with $\mathfrak{X}^*\in S^*$, $\mathfrak{X}_T\in T_h$ and $\mathfrak{X}_0$ orthogonal to $S^*$ and to $T_h$.
The set of all $\mathfrak{X}$ satisfying
\begin{align*}
    &(i)\; \mathfrak{X}^* \in \Pi\\
    &(ii)\; |\psc{\mathfrak{X}_T}{\mathfrak{Y}_j^i}| \leq H^{-(y_j-(2y-1)/(ep))}\big(H/H(B^e)\big)^{1/(2h)}\quad (j=1,\dots,h \,;\, i=1,2)\\
    &(iii)\; ||\mathfrak{X}_0||\leq C_{12}H^{(2y-1)/(ep)}
\end{align*}
(where $||\cdot||$ is the norm associated with $\psc{\cdot}{\cdot}$) is a symmetric convex body which is the product of
three symmetric convex bodies from pairwise orthogonal subspaces, hence it has a volume
\begin{align*}
    \Delta^m H(B^e)\times \Big(\prod_{j=1}^h\big(2H^{-(y_j-(2y-1)/(ep))}&\big(H/H(B^e)\big)^{1/(2h)}\big)^2\Big)\times\\
    &\times\big(C_{12}H^{(2y-1)/(ep)}\big)^{pe-2h}V(pe-2h),
\end{align*}
where $V(l)$ denotes the volume of the unit ball in $E^l$.
Finally its volume is $\Delta^m4^hV(pe-2h)C_{12}^{ep-2h}> 2^{pn}\Delta^n$ if $C_{12}$ is large enough. Therefore,
by Minkowski's Theorem, there is an $\mathfrak{X}\in\Lambda\setminus\{0\}$ in this set. One may choose $W$ in
$\mathcal{O}_{K'}^n$ such that
\[\mathfrak{X}=\rho(W).\]

\subsection{Properties of $W$}
\label{subsection Proporties of W}

 One has to establish two properties for $W$ (inequalities \eqref{going_down majoration psc(Y_i,W)} and
 \eqref{going_down majoration ||V_j||} below) in order to show that $B^{e-1}$ defined by \eqref{going_down def B^(e-1)}
 has all the required properties. More precisely, $|\psc{W}{Y_j}|$ and $||V_j||$ have to be controlled (where $V_j$ is
 the orthogonal projection of $W^{(j)}$ on $\Vect{Z_1^{(j)},\dots,Z_m^{(j)}}^{\perp}$) because these quantities will
 appear directly in the estimate of the height $H(B^{e-1})$ of $B^{e-1}$.

For $1\leq j\leq h$ one has
\begin{align}
|\psc{Y_j}{W}| &= |\psc{\mathfrak{Y}_j^1}{\rho(W)}+i\psc{\mathfrak{Y}_j^2}{\rho(W)}| = |\psc{\mathfrak{Y}_j^1}{\mathfrak{X}}+i\psc{\mathfrak{Y}_j^2}{\mathfrak{X}}| \nonumber \\
&\leq |\psc{\mathfrak{Y}_j^1}{\mathfrak{X}_T}|+|\psc{\mathfrak{Y}_j^2}{\mathfrak{X}_T}| \nonumber\\
&\leq 2H^{-(y_j-(2y-1)/(ep))}\big(H/H(B^e)\big)^{1/(2h)},
\label{going_down majoration psc(Y_i,W)}
\end{align}
by $(ii)$.\\
Also notice that $(ii)$ and $(iii)$ together imply
\[
    ||\mathfrak{X}-\mathfrak{X}^*|| \leq C_{14}H^{(2y-1)/(ep)},
\]
because by assumption $y_i \geq 1/(2h)$ (which implies that $H^{-(y_j-(2y-1)/(ep))}\times \big(H/H(B^e)\big)^{1/(2h)}$
is bounded from above by $H^{(2y-1)/(ep)}$).\\
For $1\leq j\leq p$, write $W^{(j)} = U_j+V_j$ with $U_j\in\Vect{Z_1^{(j)},\dots,Z_m^{(j)}}$ and $V_j$ orthogonal
to $\Vect{Z_1^{(j)},\dots,Z_m^{(j)}}$ (here Schmidt's arguments do not work exactly as he says, and this is what
motivates us to introduce $V_j$, $j=1,\dots,p$. See Remarks \ref{Remarque preuve Schmidt} and \ref{Remarque2 preuve Schmidt} for more details).\\
Now if $\sigma_j$ is real, then $W^{(j)},Z_1^{(j)},\dots,Z_m^{(j)} \in \RR^n$, and this implies that $U_j, V_j\in\RR^n$.
Set $\mathfrak{W}_j:= (0,\dots,\underbrace{V_j}_{j\rm{-th\;block}},\dots,0)$. Then $\mathfrak{W}_j$ is orthogonal
to $S^*$, thus to $\mathfrak{X}^*$, and to $\mathfrak{X}-\mathfrak{W}_j$ (by definition of $V_j$ and $U_j$). From
this one can deduce that
$||\mathfrak{W}_j||^2 = |\psc{\mathfrak{W}_j}{\mathfrak{X}}| = |\psc{\mathfrak{W}_j}{\mathfrak{X}-\mathfrak{X}^*}| \leq ||\mathfrak{W}_j||\times C_{14}H^{(2y-1)/(ep)}$.
Thus
\[
    ||V_j||=||\mathfrak{W}_j|| \leq C_{14}H^{(2y-1)/(ep)}.
\]
If $\sigma_j$ and $\sigma_{j+1}$ are complex conjugate, then
\[
    \mathfrak{W}_j^1:=(0,\dots,0,\underbrace{V_j^{[1]},V_j^{[2]}}_{\textrm{blocks $j$ and $j+1$}},0,\dots,0)\;\textrm{and}\; \mathfrak{W}_j^2:=(0,\dots,0,\underbrace{-V_j^{[2]},V_j^{[1]}}_{\textrm{blocks $j$ and $j+1$}},0,\dots,0)
\]
are orthogonal to $S^*$, and in particular to $\mathfrak{X}^*$. Then, one has
\begin{align*}
    &\Big|\psc{W^{[j]}}{V_j^{[1]}}+\psc{W^{[j+1]}}{V_j^{[2]}}\Big| = |\psc{\mathfrak{X}}{\mathfrak{W}_j^1}| = |\psc{\mathfrak{X}-\mathfrak{X}^*}{\mathfrak{W}_j^1}| \leq C_{14}H^{(2y-1)/(ep)}||\mathfrak{W}_j^1||,\\
    &\Big|\psc{W^{[j]}}{V_j^{[2]}}-\psc{W^{[j+1]}}{V_j^{[1]}}\Big| = |\psc{\mathfrak{X}}{\mathfrak{W}_j^2}| = |\psc{\mathfrak{X}-\mathfrak{X}^*}{\mathfrak{W}_j^2}| \leq C_{14}H^{(2y-1)/(ep)}||\mathfrak{W}_j^2||,
\end{align*}
and since
\[
    \Big|\psc{W^{(j)}}{V_j}\Big|= \Big|\psc{W^{[j]}}{V_j^{[1]}}+\psc{W^{[j+1]}}{V_j^{[2]}} +i\big(\psc{W^{[j]}}{V_j^{[2]}}-\psc{W^{[j+1]}}{V_j^{[1]}}\big)\Big|,
\]
one may conclude that
\[
    ||V_j||^2=\Big|\psc{W^{(j)}}{V_j}\Big| \leq 2C_{14}H^{(2y-1)/(ep)}||V_j||,
\]
and therefore
\begin{equation}
  \label{going_down majoration ||V_j||}
  ||V_j||\leq 2C_{14}H^{(2y-1)/(ep)}.
\end{equation}
This inequality is satisfied for $j=1,\dots,p$.\\\\

\subsection{Definition and properties of $B^{e-1}$}
\label{subsection Def et Prop B^e-1}

As announced, let $B^{e-1}$ be defined by \eqref{going_down def B^(e-1)}. Now one has to show that $H(B^{e-1})\leq H'$.
We follow Schmidt's arguments.
Let $\mathfrak{a}$ be the fractional ideal of $K'$ generated by the $\binom{n}{m}$
determinants of all $m\times m$-submatrices of the matrix with row vectors $Z_1,\dots,Z_m$, and let $\mathfrak{b}$ be
the ideal of $K'$ generated by the $\binom{n}{m+1}$ determinants of all $(m+1)\times (m+1)$-submatrices of the matrix
with row vectors $W,Z_1,\dots,Z_m$.
Since all coordinates of $W$ are in $\mathcal{O}_{K'}$ one has
$\mathfrak{b}\subset\mathfrak{a}$ (for one can use for every $(m+1)\times (m+1)$-determinant the Laplace expansion
along the first row). Thus $N(\mathfrak{b})\geq N(\mathfrak{a})$. Moreover,
$H(B^e)=H(B^{e,\perp}) = N(\mathfrak{a})^{-1}\prod_{i=1}^p ||Z_1^{(i)}\wedge\dots\wedge Z_m^{(i)}||$. Then, one has
\begin{align}
    H(B^{e-1}) &= H((B^{e-1})^{\perp}) = N(\mathfrak{b})^{-1}\prod_{i=1}^p ||W^{(i)}\wedge Z_1^{(i)}\wedge\dots\wedge Z_m^{(i)}||\nonumber\\
    & = N(\mathfrak{b})^{-1}\prod_{i=1}^p ||V_i||\;||Z_1^{(i)}\wedge\dots\wedge Z_m^{(i)}||\leq H(B^e) \prod_{i=1}^p ||V_i||\nonumber\\
    &\leq C_5H(B^e)H^{(2y-1)/e} \leq C_5H^{(e+2y-1)/e}=H',
    \label{going_down majoration H(B^(e-1)) en fonciton de H(B^e)}
\end{align}
  for a sufficiently large $C_5$.\\

Now one has to bound $H(B^{e-1})\omega_i^2(A^d,B^{e-1})$ from above (for $i=1,\dots,h$). Since
$\psc{Y_i}{Z_j}=0$ ($j=1,\dots,m$) and $||Y_i||=1$, one has
\begin{align*}
    &\dist{Y_i}{(B^{e-1})^{\perp}}^2 = \\
    & = ||Y_i\wedge W\wedge Z_1\wedge \dots\wedge Z_m||^2||Y_i||^{-2}||W\wedge Z_1\wedge\dots\wedge Z_m||^{-2}\\
    & = \Big(||Y_i||^2||W\wedge Z_1\wedge\dots\wedge Z_m||^2-|\psc{Y_i}{W}|^2||Z_1\wedge\dots\wedge Z_m||^2\Big)||W\wedge Z_1\wedge\dots\wedge Z_m||^{-2}\\
    & = 1 - |\psc{Y_i}{W}|^2||Z_1\wedge\dots\wedge Z_m||^2||W\wedge Z_1\wedge\dots\wedge Z_m||^{-2},
\end{align*}
the first equality is obtained by the special case of Proposition \ref{Prop lien mu et Det generalise} and the second one is obtained by using
Laplace's expansion twice (first for $D^2(Y_i,W,Z_1,\dots,Z_m)$ defined at the beginning of
§\ref{subsection multilinear height distance} along the first column, then for the second non-zero determinant obtained
along the first row).\\
Using Eq.~\eqref{Eq omega de l orthogonal} one finds
\[
    \dist{Y_i}{B^{e-1}} = |\psc{Y_i}{W}|\;||Z_1\wedge\dots\wedge Z_m||\;||W\wedge Z_1\wedge\dots\wedge Z_m||^{-1},
\]
hence
\[
    ||W\wedge Z_1\wedge\dots\wedge Z_m||^2\dist{Y_i}{B^{e-1}}^2 = |\psc{Y_i}{W}|^2||Z_1\wedge\dots\wedge Z_m||^2.
\]
Finally
\begin{align*}
    &H(B^{e-1})\dist{Y_i}{B^{e-1}}^2\\
    &= N(\mathfrak{b})^{-1}\Big(\prod_{j=3}^p ||W^{(j)}\wedge Z_1^{(j)}\wedge\dots\wedge Z_m^{(j)}||\Big)
    ||W\wedge Z_1\wedge\dots\wedge Z_m||^2\dist{Y_i}{B^{e-1}}^2\\
    &=N(\mathfrak{b})^{-1}\Big(\prod_{j=3}^p ||W^{(j)}\wedge Z_1^{(j)}\wedge\dots\wedge Z_m^{(j)}||\Big)
    ||Z_1\wedge\dots\wedge Z_m||^2|\psc{Y_i}{W}|^2\\
    &\leq C_{15}N(\mathfrak{a})^{-1}\Big(\prod_{j=1}^p ||Z_1^{(j)}\wedge\dots\wedge Z_m^{(j)}||\Big)
    \Big(\prod_{j=3}^p||V_j||\Big)H^{-2(y_i-(2y-1)/(ep))}\big(H/H(B^e)\big)^{1/h}\\
    &\leq C_{16}H(B^e)H^{(p-2)(2y-1)/(ep)-2[y_i-(2y-1)/(ep)]}H/H(B^e)\\
    &\leq C_{16}H^{(e+2y-1)/e-2y_i}.
\end{align*}
The first inequality follows from \eqref{going_down majoration psc(Y_i,W)} and the second one
from \eqref{going_down majoration ||V_j||}. Then, there is a vector $R_i\in B^{e-1}\setminus\{0\}$ having
\[
    H(B^{e-1})\dist{Y_i}{R_i}^2 \leq C_{16}H^{-(2y_i'-1)(2y+e-1)/e},
\]
since by definition $y'_i:= (y_ie)/(qy+e-1)$.\\
By assumption, one has
\[
    H(B^e)\dist{X_i}{Y_i}^2 \leq c^2H^{-2y_i+1},
\]
thus, by \eqref{going_down majoration H(B^(e-1)) en fonciton de H(B^e)} :
\[
    H(B^{e-1})\dist{X_i}{Y_i}^2 \leq c^2C_5H^{-2y_i+1+(2y-1)/e}=c^2C_5H^{-(2y_i'-1)(2y+e-1)/e}.
\]

These inequalities and the triangle inequality provide
\[
    H(B^{e-1})\dist{X_i}{R_i}^2 \leq C_{17}c^2H^{-2y_i+1+(2y-1)/e}=C_{17}c^2H^{-(2y_i'-1)(2y+e-1)/e}\quad(i=1,\dots,h),
\]
and so \cite[Theorem $7$]{Schmidt67} (with $\delta = 1$ for instance, since $(X_j)_j$ is an orthonormal family)
yields
\[
    H(B^{e-1})\omega_i^2(A^d,B^{e-1}) \leq c^2C_6H^{-(2y_i'-1)(2y+e-1)/e}\quad(i=1,\dots,h),
\]
for $C_6$ large enough. This completes the first part of the proof.\\\\

\subsection{Proof of the second part of the theorem}
\label{section de la remark corrective}

Suppose now that \eqref{going_down 3ème equation} holds. We follow Schmidt's proof and define $H_1\geq 1$ by
$H'=C_9HH_1^{2h/e}$ ($C_9$ will be specified later). For the construction of $W$, replace $(i)$, $(ii)$, $(iii)$ by
\begin{align*}
    &(i')\; \mathfrak{X}^* \in \Pi\\
    &(ii')\; |\psc{\mathfrak{X}_T}{\mathfrak{Y}_j^i}| \leq H_1^{-(1-2h/(ep))}\quad (j=1,\dots,h \,;\, i=1,2)\\
    &(iii')\; ||\mathfrak{X}_0||\leq C_{18}H_1^{2h/(ep)}.
\end{align*}
These equations define a symmetric convex body of volume
\begin{align*}
    &\Delta^mH(B^e)\times (2H_1^{-(1-2h/(ep))})^{2h}\times (C_{18}H_1^{2h/(ep)})^{ep-2h}V(ep-2h)\\
    & = \Delta^mH(B^e)4^hC_{18}^{ep-2h}V(ep-2h)\\
    & > 2^{pn}\Delta^n,
\end{align*}
for $C_{18}$ large enough. By Minkowski's Theorem there is an $\mathfrak{X}\in\Lambda\setminus\{0\}$ in this set.
Let $W$ be in $\mathcal{O}_{K'}^n$ such that $\mathfrak{X}=\rho(W)$.\\
Equation \eqref{going_down majoration psc(Y_i,W)} is replaced by
\[
    |\psc{W}{Y_j}|\leq 2H_1^{-(1-2h/(ep))}.
\]
One also has
\[
    ||\mathfrak{X}-\mathfrak{X}^*||\leq C_{19}H_1^{2h/(ep)},
\]
thus \eqref{going_down majoration ||V_j||} becomes
\[
    ||V_j|| \leq 2C_{19}H_1^{2h/(ep)}.
\]
Computing $H(B^{e-1})$ one finds
\[
    H(B^{e-1})\leq H(B^e)\prod_{j=1}^p||V_j|| \leq CH(B^e)H_1^{2h/e}.
\]
Now set $C_9:= C$, which implies $H(B^{e-1})\leq H'$. With these new estimates one finds
\begin{align*}
    H(B^{e-1})\dist{Y_i}{B^{e-1}}^2 &\leq C_{20}H(B^e)\times H_1^{(p-2)2h/(ep)}\times H_1^{-2(1-2h/(ep))}\\
    &\leq C_{21} HH_1^{2h/e-2} = C_{21} HH_1^{2h/e(1-e/h)} \\
    &= C_{10}H^{e/h}H'^{1-e/h}=C_{10}H^{2y'_0}H'^{-(2y'_0-1)}.
\end{align*}
Now, note that \eqref{going_down 3ème equation} implies that $Y_i = \pm X_i$ (for $i=1,\dots,h$) and
so $Y_1,\dots,Y_h$ form an orthonormal subset of $A^d$; \cite[Theorem $7$]{Schmidt67} (with $\delta = 1$)
yields the expected result.\\

\end{proof}

\begin{Rem}
\label{Remarque preuve Schmidt}
We give here more details in case $(b)$ about the reasons which lead us to introduce subspaces $V_j$ and the decomposition of
vectors $W^{(j)}$ (see Subsection \ref{subsection Proporties of W}). In his paper \cite[p.~455]{Schmidt67},
Schmidt writes (with his notation for the orthogonal complement and  for $Z_j\in K^n$) $W=U+V$, where
$U\in B^{e\perp} = \Vect{Z_1,\dots, Z_m}$ and $V\in B^e$, in order to have the decomposition
$W^{(j)} = U^{(j)} + V^{(j)}$ with $U^{(j)}\in (B^e)^{(j)\perp} = \Vect{Z_1^{(j)},\dots, Z_m^{(j)}}$ and
$V^{(j)}\in (B^e)^{(j)}$.\\
Here a problem arises : in the complex case, which orthogonal complement does the symbol $\perp$ denote ? With our
notation, suppose that one considers $(B^e)^{\phi,\perp}$. Then, it is not always true that we have
$G^n = B^e \oplus (B^e)^{\phi,\perp}$ and so the decomposition $W=U+V$ may not be considered. (This is the
first but not the only one problematic point: for example with this definition it seems that we could not obtain
 inequality $(15)$ of Schmidt, because just before in his text the notation $W^{(j)}R^{(j)}$ would refer to $\phi(W^{(j)},R^{(j)})$ and not to the inner product
$\psc{W^{(j)}}{R^{(j)}}$).\\
Suppose now that $(B^e)^{\perp}$ denotes the orthogonal complement for the inner product $\psc{\cdot}{\cdot}$. Here,
the decomposition $W=U+V$ may be considered. However it rises another problem : $(B^e)^{\perp}$ is no longer defined
over $K$, but over $K'$, and we could not apply the embedding $\sigma_j : K \rightarrow \CC$ to $U$.
In fact we also have a problem to define $V^{(j)}$ although $B^e$ is defined over $K$, because we do not have
necessarily $V\in K^n$ (the decomposition $(0,1) = \lambda(1,\alpha) + \lambda(-1,(\overline{\alpha})^{-1})$ with
$\lambda = 1/(\alpha + (\overline{\alpha})^{-1})$ and $\alpha\in K\setminus{K'}$ provides a simple
counter-example in dimension two). Furthermore, even if $K = K'$ and that $U^{(j)}$ and $V^{(j)}$ could be
considered, it is not true that it implies $\psc{U^{(j)}}{V^{(j)}} = 0$, because of the complex conjugation in the
inner product. This brings problems to apply Schmidt's arguments.\\

\end{Rem}

\begin{Rem}
\label{Remarque2 preuve Schmidt}
In case $(a)$, it makes sense to considerer the decomposition $W^{(j)} = U^{(j)} + V^{(j)}$ of Schmidt, but if $K$ is
not totally real and if $\sigma_j$ is a non-real embedding, with Schmidt's argument it seems that one can only obtain
\[
    |\phi(V^{(j)},V^{(j)})| \leq 2C_{14}H^{(y-1)/(ep)} ||V^{(j)}||,
\]
instead of inequality $(15)$ of \cite{Schmidt67} (which corresponds to inequality \eqref{going_down majoration ||V_j||}
in this paper).\\
To solve this problem, it suffices to apply the argument used for the complex case. It is simpler in case $(a)$ because
we have $K'=K$, and $B^e$ is the subspace defined by Schmidt in \cite{schanuel1964heights} (defined over $K$).
If we replace the decomposition $U^{(j)}=W^{(j)}+U^{(j)}$ with the decomposition $W^{(j)} = U_j+V_j$ if
$\sigma_j$ is not real, with $U_j\in\Vect{Z_1^{(j)},\dots,Z_m^{(j)}}$ and $V_j$ orthogonal
to $\Vect{Z_1^{(j)},\dots,Z_m^{(j)}}$, as for \eqref{going_down majoration ||V_j||} we obtain
\[
    ||V_j|| \leq 2C_{14}H^{(y-1)/(ep)}.
\]
This inequality is slightly different from \eqref{going_down majoration ||V_j||} because in case $(a)$ the
symmetric convex body defined by $(i)$, $(ii)$ and $(iii)$ at the beginning of the proof is not the same.
Now, note that in Section \ref{subsection Def et Prop B^e-1} we use inequalities \eqref{going_down majoration psc(Y_i,W)} and
\eqref{going_down majoration ||V_j||} but not directly the fact that $K\subset\RR$ or not, so working with $V_j$ rather
than $V^{(j)}$ we may follow Schmidt's arguments to complete the proof of case $(a)$.

\end{Rem}

\section{Exponents of Diophantine Approximation}
\label{section application exponents}

Let $n\geq 1$ and $K\subset\CC$ be a number field. We recall that one distinguishes between two cases $(a)$ and $(b)$.\\
In case $(a)$, $K$ is real, $G^{n+1}$ denotes Euclidean space $\RR^{n+1}$, and $q=1$.\\
In case $(b)$, $K$ is complex non real, $G^{n+1}$ denotes unitary $\CC^{n+1}$, and $q=2$.

\begin{Def}
Let $\bu\in G^{n+1}\setminus\{0\}$. For each $j=0,\dots,n-1$, we denote by $\wj{j}{K}{\bu}$ (resp. $\wjhat{j}{K}{\bu}$)
the supremum of all real numbers $\omega$ such that, for arbitrarily large values of $Q$ (resp. for all sufficiently
large values of $Q$), there exists a vector subspace $S$ of $G^{n+1}$, defined over $K$, of dimension $j+1$, with
\[
    H(S)\leq Q \quad \textrm{and}\quad H(S)\omega_1^q(\bu,S)\leq Q^{-\omega}.
\]
\end{Def}

Laurent introduces this family of exponents in \cite{laurent703146transfer} in the case $K=\QQ$. He gives a series of
inequalities (which may be proved using Schmidt's results \cite{Schmidt67}) relating these exponents (cf for
example \cite[Theorem 2.2]{Roy_octobre}, and compare to Theorem \ref{thm encadrement omega_j} below), and gives also
a description of the full spectrum of the $2n$ exponents
$(\wj{0}{\QQ}{\bu},\dots,\wj{n-1}{\QQ}{\bu},\dots,\wjhat{0}{\QQ}{\bu},\dots,\wjhat{n-1}{\QQ}{\bu})$ for $n=2$.
In \cite{Roy_octobre}, Roy gives a description of the full spectrum of the $n$ exponents
$(\wj{0}{\QQ}{\bu},\dots,\wj{n-1}{\QQ}{\bu})$ for all $n\geq 1$, proving that the inequalities of Laurent describe
completely this spectrum (cf Theorem 2.3 of \cite{Roy_octobre}).\\

For an arbitrary number field $K$, the aforementioned inequalities can be generalized: this is
Theorem \ref{thm encadrement omega_j}, which is the main application of this section.

\begin{Thm}
\label{thm 13}
Let $\bu \in G^{n+1}\setminus\{0\}$ and let $j$, $0\leq j\leq n-1$. Then
\begin{align}
\wj{j}{K}{\bu} \geq \wjhat{j}{K}{\bu} \geq \frac{j+1}{n-j}\quad (0\leq j \leq n-1).
\end{align}
\end{Thm}

\begin{Dem}[of Theorem \ref{thm 13}]\newline
By \cite[Theorem $13$]{Schmidt67} in the case $d=1$, there exists a constant $c>0$ which depends only of $n$ and $K$
such that for all $Q\geq 1$, there is a subspace $S$, defined over K and of dimension $j+1$, satisfying
\[
    H(S) \leq Q \quad \textrm{and} \quad H(S)\omega_1^q(\bu,S) \leq cQ^{-(j+1)/(n-j)}.
\]
Theorem \ref{thm 13} follows immediately. Note that in the proof of Theorem $13$, Schmidt uses his going-down Theorem.

\end{Dem}

\begin{Thm}
\label{thm encadrement omega_j}
Let $n\in \NN^*$ and $\bu \in G^{n+1}$. Then we have $\wj{0}{K}{\bu}\geq \frac{1}{n}$ and
\begin{align}
\label{inegalites omega_j,K(u)}
\frac{j\wj{j}{K}{\bu}}{\wj{j}{K}{\bu}+j+1} \leq \wj{j-1}{K}{\bu} \leq \frac{(n-j)\wj{j}{K}{\bu}-1}{n-j+1}\quad (1\leq j\leq n-1),
\end{align}
with the convention that the left-most ratio is equal to $j$ if $\wj{j}{K}{\bu} = \infty$.
\end{Thm}

The right inequality in \eqref{inegalites omega_j,K(u)} follows from Schmidt's going-up
Theorem \cite[Theorem 9]{Schmidt67}, while the left inequality follows from the going-down Theorem. In his
paper \cite{laurent703146transfer}, Laurent introduces the exponents $\wj{j}{\QQ}{\bu}$ and notes that each
inequalities in \eqref{inegalites omega_j,K(u)} is best possible (because they allow to find Khinchine's
transference inequalities, which are best possible). For $\bu$ with $\QQ-$linearly independent coordinates
and $K=\QQ$, Laurent gives an independent proof of the right inequality of \eqref{inegalites omega_j,K(u)}
in \cite{laurent703146transfer}, and both inequalities are proved by Bugeaud and Laurent in \cite{bugeaud2010transfer}.
Roy proves that for $K=\QQ$, the inequalities of Theorem \ref{thm encadrement omega_j} describe the set of all
possible values of the $n$-tuples $(\wj{0}{\QQ}{\bu},\dots,\wj{n-1}{\QQ}{\bu})$ (cf \cite[Theorem 2.3]{Roy_octobre}).
It would be of interest to know if his theorem remains true or not for an arbitrary number field $K$.\\

\begin{Dem}

We show first the left inequality in \eqref{inegalites omega_j,K(u)}.

Set $A^1 := \textrm{Span}(\bu)$, the line spanned by $\bu$ in $G^{n+1}$.
Let $j$, $1\leq j \leq n-1$. If $\wj{j}{K}{\bu} = \infty$ let $y_1 > q^{-1}$. Otherwise let $y_1= q^{-1}(\wj{j}{K}{\bu}+1) - \varepsilon$, where $\varepsilon > 0$ is
very small; then $y_1>q^{-1}$ because Theorem~\ref{thm 13} yields $\wj{j}{K}{\bu} > 0$. Then, there exist
arbitrarily large values of $Q$ for which there exists a subspace $S$ of $G^{n+1}$ of dimension $j+1$ ($1\leq j \leq n-1$), such that
\begin{align}
\label{equation 1}
    H(S)\leq Q \quad \textrm{and}\quad H(S)\omega_1^q(A^1,S) \leq Q^{-(qy_1-1)}.
\end{align}

Set $y'_1 := y_1(j+1)/(qy_1+j)$. The conditions of the going-down Theorem are fulfilled with $A^1$ and $B^e := S$
(one has $d = i = h = 1$, $e = j+1$ and we choose $c=1$). This gives a subspace $S'\subset S$ of height $H(S')\leq Q'$
with $Q':=C_5Q^{\frac{qy_1+j}{j+1}}$, of dimension $j$, defined over $K$, having
\begin{align}
\label{equation 2}
    H(S')\leq Q'\quad \textrm{and}\quad H(S')\omega_1^q(\bu,S') \leq C_7 Q'^{-(qy'_1-1)}.
\end{align}

Since Equation \eqref{equation 1} holds for arbitrarily larges values of $Q$, \eqref{equation 2} holds also for
arbitrarily large values of $Q'$, and one deduces that
\[
    \wj{j-1}{K}{\bu} \geq qy'_1-1.
\]
Letting $y_1$ tends to $q^{-1}(\wj{j}{K}{\bu} + 1)$ if $\wj{j}{K}{\bu}<\infty$, and to $+\infty$ otherwise, we deduce the left side of the inequality in \eqref{inegalites omega_j,K(u)}.\\

Only a scheme of proof is given for the right inequality.\\
Set again $A^1 := \textrm{Span}(\bu)$ and let $j$, $1\leq j \leq n-1$. We can suppose that for any subspace $S$ of
dimension $j+1$, defined over $K$, we have $\omega_1^q(A^1,S) > 0$ (otherwise it means that $\wj{j}{K}{\bu} = \infty$,
and the right inequality in \eqref{inegalites omega_j,K(u)} is obvious). Let $y_1$, $0<qy_1<\wj{j-1}{K}{\bu}$. There
exist arbitrarily large values $Q\geq 1$ and corresponding subspaces $S$ of dimension $j$ defined
over $K$, having
\begin{align}
\label{equation 3}
    H(S)\leq Q \quad \textrm{and}\quad H(S)\omega_1^q(A^1,S) \leq Q^{-qy_1}.
\end{align}
Now, we use Schmidt's going-up Theorem (\cite[Theorem 8]{Schmidt67}) with parameters $n+1$, $x_1 = 1/q$,
$y_1$, $d=t=i=1$, $c=1$, and $B^e := S$ (note that in our context we have $\psi_1(A^1,B^e) = \omega_1(A^1,B^e)$ for all
$B^e$ of dimension $e < n$). This gives us a subspace $S' \supseteq S$ of dimension $j+1$, defined over $K$, having
\begin{align}
\label{equation 4}
    H(S') \leq Q' \quad \textrm{and} \quad 0 < H(S')^{(n-j+1)/(n-j)}\omega_1^q(A^1,S') \leq C_4 Q'^{-qy_1(n-j+1)/(n-j)},
\end{align}
where  $Q':= C_3Q^{(n-j)/(n-j+1)} $ and $C_3,C_4>0$ depend of $n, K$ and $y_1$ only.\\
Now, since Equation \eqref{equation 3} holds for arbitrarily large values $Q$, so does \eqref{equation 4} for
arbitrarily large values $Q'$. This implies that there are infinitely many subspaces $S'$ which satisfy
\eqref{equation 4}. In particular, $H(S')$
is not bounded from above.\\
To conclude, it suffices to remark that \eqref{equation 4} implies
\[
    0 < H(S')^{(n-j+1)/(n-j)}\omega_1^q(A^1,S') \leq C_4 H(S')^{-qy_1(n-j+1)/(n-j)}.
\]
Then multiplying each side of the inequality by $H(S')^{1-(n-j+1)/(n-j)}$ and writing
$Q''=H(S')$, we find that for arbitrarily large values of $Q''$ there exists $S'$ of dimension $j+1$ defined over
$K$ such that
\begin{align}
\label{equation 5}
    H(S') \leq Q'' \quad \textrm{and} \quad H(S')\omega_1^q(A^1,S') \leq C_4 Q''^{-(qy_1+1)(n-j+1)/(n-j) + 1}.
\end{align}
This shows that
\[
    \wj{j}{K}{\bu} \geq (qy_1+1)(n-j+1)/(n-j) - 1,
\]
which gives
\[
    qy_1 \leq \frac{(n-j)\wj{j}{K}{\bu} - 1}{n-j+1},
\]
and letting $qy_1$ tend to $\wj{j-1}{K}{\bu}$ we find the desired result.

\end{Dem}

\bibliographystyle{abbrv}  


\Addresses

\end{document}